\magnification=1200
\centerline
{\bf Cohomological invariants for $G$-Galois algebras and self-dual normal bases}
\bigskip
\centerline
{E. Bayer-Fluckiger and R. Parimala}
\bigskip
\bigskip

{\bf Abstract.} We define degree two cohomological invariants for $G$-Galois algebras over fields of characteristic not 2, and use them to give necessary conditions
for the existence of a self--dual normal basis. In some cases (for instance, when the field has cohomological dimension $\le 2$) we show that these conditions are
also sufficient.

\bigskip
{\bf Introduction}
\bigskip
Let $k$ be a field of characteristic $\not = 2$, and let $L$ be a finite degree Galois extension of $k$. 
Let $G = {\rm Gal}(L/k)$. The {\it trace form} of $L/k$ is by definition the quadratic form
$q_L : L \times L \to k$ defined by
$q_L(x,y) = {\rm Tr}_{L/k}(xy)$.  Note that $q_L$ is a {\it $G$-quadratic form}, in other
words we have $q_L(gx,gy) = q_L(x,y)$ for all $x,y \in L$.
A normal basis $(gx)_{g \in G}$ of $L$ over $k$ is said to be 
 {\it self-dual} 
if $q_L(gx,gx)= 1$ and $q_L(gx,hx) = 0$ if $g \not = h$. It is natural to ask which extensions have a self-dual
normal basis. This question is investigated
in several papers (see for instance [BL 90], [BSe 94],  [BPS 13]). It is necessary to work in a more general context than the one of Galois extensions, namely that of
$G$-{\it Galois algebras} (see for instance [BSe 94], \S 1);  one advantage being that this category is stable by base change of the ground field;
the notion of a self-dual normal basis is defined in the same way. 

\medskip
If $k$ is a global field, then the {\it Hasse principle} holds  : a $G$-Galois algebra has a self-dual normal basis over $k$ if
and only if such a basis exists everywhere locally (see [BPS 13]). The present paper completes this result by giving necessary
and sufficient conditions for the existence of a self-dual normal basis when $k$ is a local field (cf. \S 7). The conditions are given in terms
of cohomological invariants  defined over the ground field $k$ constructed in \S 3 and \S 4.

\medskip
For an arbitrary ground field $k$, we start with the $H^1$-invariants defined in [BSe 94], \S 2. Recall from [BSe 94] that the vanishing of these invariants is a necessary condition for the
existence of a self-dual normal basis; it is also sufficient in the case of fields of cohomological dimension 1 (see [BSe 94], Corollary 2.2.2 and
Proposition 2.2.4). 

\medskip Let $k[G]$ be the group algebra of $G$ over $k$, and let $J$ be its radical; the quotient $k[G]^s = k[G]/J$ is a semisimple $k$-algebra. Let $\sigma : k[G] \to k[G]$ be the $k$-linear involution sending $g$ to $g^{-1}$; it induces an involution $\sigma^s : k[G]^s \to k[G]^s$. The algebra $k[G]^s$ splits as a product of simple algebras.  If $A$ is a 
$\sigma^s$-stable simple algebra which is a 
factor of $k[G]^s$, we denote by $\sigma_A$ the
restriction of $\sigma^s$ to $A$, and by $E_A$ the subfield of the center of $A$ fixed by $\sigma_A$. We say that $A$ is  {\it orthogonal} if $\sigma_A$ is
the identity on the center of $A$, and if over a separable closure of $k$  it is induced by a symmetric form, and  {\it unitary} if $\sigma_A$ is not the identity
on the center of $A$ (see 1.3 for details). 

\medskip Let $L$ be a $G$-Galois algebra over $k$, and let us assume that its $H^1$-invariants are trivial. 
We then define, for every orthogonal or unitary  $A$ as above, 
cohomology classes in $H^2(k, {\bf Z}/2{\bf Z})$, denoted by $c_A(L)$ in the orthogonal case and by $d_A(L)$ in the unitary case (see
\S 3 and \S 4). They are invariants of the $G$-Galois algebra $L$. They also provide necessary conditions for the existence of a self-dual normal
basis (this involves restriction to certain finite degree extensions of $k$, namely, the extensions  $E_A/k$; see Propositions 3.5 and 4.7 for precise statements).
If moreover $k$ has cohomological dimension $\le 2$, then these conditions are also sufficient  
(Theorem 5.3.). Finally, if $k$ is a local field, then the conditions can be expressed in terms of the invariants $c_A(L)$ and $d_A(L)$, without passing
to finite degree extensions (Theorem 7.1).  Section 8 applies the results of \S 7 and the Hasse principle of [BSP 13] to give necessary
and sufficient conditions for the existence of a self-dual normal basis when $k$ is a global field (Theorem 8.1).  

\medskip
Section 6 deals with the case of cyclic groups of order a power of $2$ over arbitrary fields. We show that at most one of the unitary components $A$ gives rise to
a non-trivial invariant $d_A(L)$ (Proposition 6.4 {\rm (i)}), and that this invariant provides a necessary and sufficient condition for the existence of a self-dual normal basis
(Corollary 6.5). 


\bigskip
\bigskip
{\bf \S 1. Definitions, notation and basic facts}
\bigskip
{\bf 1.1. Galois cohomology}
\medskip
We use standard notation in Galois cohomology. If $K$ is a field, we denote by $K_s$ a separable closure of $K$, and by $\Gamma_K$ the Galois group ${\rm Gal}(K_s/K)$. 
For
any discrete $\Gamma_K$-module $C$, set $H^i(K,C) = H^i(\Gamma_K,C)$. If $\Gamma$ is a finite or profinite group, set $H^i(\Gamma) = H^i(\Gamma,{\bf Z}/2{\bf Z})$.
If $U$ is a $K$-group scheme, we denote by
 $H^1(K,U)$ the pointed set $H^1(\Gamma_K,U(K_s))$. 
 \bigskip
{\bf 1.2. Algebras with involution and unitary groups}
\medskip
Let $K$ be a field of characteristic $\not = 2$, and 
let $R$ be a finite dimensional algebra over $K$. An {\it involution} of $R$ is a
$K$-linear anti-automorphism
$\sigma : R \to R$ such that $\sigma^2$ is the identity. 


\medskip Let us denote by ${\rm Comm}_K$ the category of commutative $K$-algebras, and by ${\rm Group}$ the category of groups. If $(R,\sigma)$ is an algebra with involution, the functor
${\rm Comm}_K \to {\rm Group}$ given by $S \mapsto  \{x \in R \otimes_K S  \ | \ x \sigma( x )= 1 \}$ is the functor of points of a scheme over 
${\rm Spec} (K)$; we denote it by $U_{R,K}$.

\medskip
Let $h = \langle 1 \rangle$ be the rank one unit hermitian form over $(R,\sigma)$, given by $h(x,y) = x \sigma(y)$ for all $x,y \in R$. Then
$U_{R,K}$ is  the scheme of automorphisms of the hermitian form  $h$. This is a smooth, finitely
presented affine group scheme over ${\rm Spec} (K)$ (see for instance [BF 15], Appendix A).
Moreover, $H^1(K,U_{R,K})$ is in natural bijection with the set of isomorphism classes of rank one hermitian forms
over $(R,\sigma)$ that become isomorphic to $h$ over $K_s$ (see [Se 64], chap. III, \S 1). 

\medskip
If $F$  is a subfield of $K$, then $U_{R,F} = {\rm R}_{K/F}(U_{R,K})$, where ${\rm R}_{K/F}$ denotes Weil restriction of scalars relative to
the extension $K/F$.

\medskip
Let $Z$ be the center of $R$, and assume that $R$ is a simple algebra. We say that $(R,\sigma)$
is a {\it central simple algebra with involution over $K$} if the fixed field of $\sigma$ in $Z$
is equal to $K$. If $(R,\sigma)$ is central simple algebra with involution over $K$, we set $U_R = U_{R,K}$.


\bigskip

{\bf 1.3. D\'evissage}

\bigskip
Let $G$ be a finite group and  let $k[G]$ be its group algebra over $k$. The {\it canonical involution} of $k[G]$ is the 
$k$-linear involution $\sigma : k[G] \to k[G]$ such that $\sigma(g) = g^{-1}$ for all $g \in G$.  Let $J$ be the radical of $k[G]$, and set $k[G]^s = k[G]/J$; it is a semisimple
$k$-algebra. Since 
$J$ is stable by $\sigma$, we obtain an involution $\sigma^s : k[G]^s \to k[G]^s$. Set $U_G = U_{k[G],k}$ and $U_G^r = U_{k[G]^s,k}$.
Let $N$ be the kernel of the natural surjection $U_G \to U^r_G$. Let us define group schemes $N_i$ by setting $N_i(S) = \{x \in N(S) \ | \ x \equiv 1 \ {\rm mod} \ J^i \otimes_k S \}$.
Then $1 = N_m \subset N_{m-1} \subset \cdots \subset N_1 = N$, where $m$ is an integer such that $J^m = 0$. Note that $J^i/J^{i+1}$ is a module over the semisimple
algebra $k[G]^s$, hence $N_i/N_{i-1}$ is isomorphic to a finite product of additive groups ${\bf G}_a$; therefore $N$ is a split unipotent group. This implies that
$H^1(k,U_{G})  = H^1(k,U^r_{G})$ (see for instance [Sa 81], Lemme 1.13). 

\medskip

The semisimple algebra $k[G]^s$ is known to be a direct product of simple algebras.
Note that $k[G]$ comes by scalar extension  from $k_0[G]$ for $k_0 = {\bf Q}$ or ${\bf F}_p$, hence the centers of the factors of $k[G]^s$ are  abelian Galois extensions of  $k$  of finite degree;
some  are stable under $\sigma^s$ (we call them $A$), and others come in pairs, interchanged by $\sigma^s$ (we call them $B$).

\medskip 
If $A$ is a $\sigma^s$-stable  simple factor of $k[G]^s$, we denote by $\sigma_A$ the restriction of $\sigma^s$ to $A$, by $F_A$ the center of $A$, and
by $E_A$ the subfield of $\sigma_A$-invariant elements of $F_A$. Note that $U_A$ is a group scheme over ${\rm Spec}(E_A)$.
Similarly, if
$B$ is the product of two simple algebras interchanged by $\sigma^s$, we denote by $E_B$ the subfield of the center of $B$ fixed by the involution;
$U_{B,E_B}$ is a group scheme over ${\rm Spec}(E_B)$.

\medskip We have $U^r_{G} \simeq \prod_{A} {\rm R}_{E_A/k}(U_A) \times \prod_B {\rm R}_{E_B/k} (U_{B,E_B})$, hence 



$$H^1(k,U^r_{G}) = \prod_{A} H^1(k, {\rm R}_{E_A/k}(U_A)) \times \prod_B  H^1(k,{\rm R}_{E_B/k}( U_{B,E_B})).$$

Note that $H^1(k,{\rm R}_{E_B/k}( U_{B,E_B})) = H^1(E_B,U_{B,E_B}) = 0$  (see for instance [KMRT 98], (29.2)), that 
$H^1(k,{\rm R}_{E_A/k}(U_A)) = H^1(E_A,U_A)$ (see for instance [O 84], 2.3), and that $H^1(k,U_{G})  = H^1(k,U^r_{G})$ (see above). Therefore we have 
$$H^1(k,U_{G}) = \prod_{A} H^1(E_A,U_A).$$




The algebras with involution $(A,\sigma_A)$ appearing in this product are of three types :

\medskip
{\rm (a)} The involution $\sigma_A : A \to A$ is not the identity on the center $F_A$ of $A$. Hence  $F_A/E_A$ is
a quadratic extension. Such an algebra with involution is called {\it unitary}; the group scheme $U_A$ is of Dynkin type A.

\medskip
{\rm (b)} The involution $\sigma_A : A \to A$ is the identity on $F_A$  (which is then equal to $E_A$), and, over a separable
closure of $E_A$, the involution is induced by a symmetric form. Such an algebra with involution is called {\it orthogonal}; the group scheme
$U_A$ is of Dynkin type B or D. 

\medskip
{\rm (c)} The involution $\sigma_A : A \to A$ is the identity on $F_A$ (which is then equal to $E_A$), and, over a separable
closure of $E_A$,  the involution is induced by a skew-symmetric form. Such an algebra with involution is called {\it symplectic}; the group scheme
$U_A$ is of Dynkin type C.

\bigskip
{\bf 1.4. $G$-quadratic forms}
\medskip
A {\it $G$-quadratic form} is a pair $(M,q)$, where $M$ is a $k[G]$-module
that is a finite dimensional $k$-vector space, 
and $q : M \times M \to k$ is a non-degenerate symmetric bilinear
form such that $$q(gx,gy) = q(x,y)$$ for all $x,y \in M$ and all $g \in G$.  We say that
two $G$-quadratic forms $(M,q)$ and $(M',q')$ are {\it isomorphic} if there exists
an isomorphism of $k[G]$-modules $f : M \to M'$ such that $q'(f(x),f(y)) = q(x,y)$ for
all $x,y \in M$. If this is the case, we write $(M,q) \simeq_G (M',q')$, or $q \simeq_G q'.$
It is well-known that $G$-quadratic forms correspond bijectively to non-degenerate hermitian forms
over $(k[G],\sigma)$  (see for
instance [BPS 13], 2.1, Example on page 441). 
The {\it unit $G$-form} is by definition the pair $(k[G],q_0)$, where
$q_0$ is the $G$-form characterized by $q(g,g) = 1$ and $q(g,h) = 0$ if $g \not = h$, for $g,h \in G$. 

\bigskip
{\bf 1.5. Trace forms of $G$-Galois algebras}
\medskip
If  $L$ is an \'etale $k$-algebra, we denote by 
$$q_L : L \times L \to k, \ \ q_L(x,y) = {\rm Tr}_{L/k}(xy),$$
its trace form. Then $q_L$ is a non-degenerate quadratic form over $k$; if moreover $L$ is a $G$-Galois algebra, then $q_L$ is a $G$-quadratic form. 

\medskip Let $L$ be a $G$-Galois algebra; then
$L$ has a self-dual normal basis over $k$ if and only
if $q_L$ is isomorphic to $q_0$ as a $G$-quadratic form. 
Let $\phi : \Gamma_k \to G$ be a continuous
homomorphism corresponding to $L$ (see for instance [BSe 94], 1.3). Recall that $\phi$ is unique up to conjugation. The composition 
$$\Gamma_k \  {\buildrel {\phi} \over \to} \  G \to U_{G}(k) \to U_{G}(k_s)$$
is a $1$-cocycle $\Gamma_k \to U_G(k_s)$.
Let $u(L)$ be its class in the cohomology set $H^1(k,U_G)$; it does not depend on the choice of $\phi$. The $G$-Galois algebra
$L$ has a self-dual normal basis over $k$ if and only if $u(L) = 0$, cf. [BSe 94], Corollaire 1.5.2.

\medskip
Recall from 1.3 that we have $$H^1(k,U_{G})   = \prod_{A} H^1(E_A,U_A).$$ Let $u_A(L)$ be the image of  $u(L)$ in $H^1(E_A,U_A)$; note that  $L$ has
a self-dual normal basis if and only if $u_A(L) = 0$ for every $A$.

\medskip Let $A$ be as above. 
Composing the injection $G \to U_{G}(k)$ with the natural map $U_{G}(k)  \to U^r_{G}(k) \to {\rm R}_{E_A/k}(U_{A})(k) = U_A(E_A)$, we obtain a homomorphism
$G \to U_A(E_A)$, denoted by $i_A$.

\medskip
Let $\phi_{A} : \Gamma_{E_{A}} \to \Gamma_k \to G$ be the composition of $\phi : \Gamma_k \to G$ with the inclusion of 
$\Gamma_{E_A}$ in $\Gamma_k$.  Composing  $\phi_{A}$ with the map
$i_A : G \to U_A(E_A)$ defined above we obtain a $1$-cocycle $\Gamma_{E_A} \to U_A(k_s)$. The class of this $1$-cocycle in $H^1(E_A,U_A)$
is equal to $u_A(L)$.

\bigskip
\bigskip
{\bf  \S 2. The $H^1$-condition}
\bigskip
Let $L$ be a $G$-Galois algebra over $k$, and let 
$\phi : \Gamma_k \to G$ be a homomorphism  corresponding to $L$.  Let $n$ be an integer $\ge 1$. 
 Then $\phi$ 
 induces a homomorphism
$$\phi^* : H^n(G) \to H^n(k,{\bf Z}/2{\bf Z}).$$  Note that $\phi^*$ is independent of the choice of $\phi$ in 
its conjugacy class (see [Se 68], chap. VII, proposition 3). For all $x \in H^n(G)$, set $x_L = \phi^*(x)$. 
\medskip
\noindent
{\bf Proposition 2.1.} {\it If $L$ has a self-dual normal basis over $k$, then
for all $x \in H^1(G)$ we have $x_L = 0$.}
\medskip
\noindent
{\bf Proof.} See [BSe 94], Corollaire 2.2.2. 
\medskip
If  ${\rm cd}_2(\Gamma_k) \le 1$, then $L$ has a self-dual normal basis over $k$ if and only
if $x_L = 0$ 
for all $x \in H^1(G)$, 
see [BSe 94],  Proposition 2.2.4.

\medskip
We say that {\it the $H^1$-condition is satisfied} if
$x_L = 0$ 
for all $x \in H^1(G)$. 
Let $G^2$ be the subgroup of $G$ generated
by the squares of elements of $G$. Note that $G/G^2$ is an elementary
abelian 2-group, and  that the $H^1$-condition means that the homomorphism
$\Gamma_k \to G \to G/G^2$ induced by $\phi$ is trivial, i.e.  $\phi(\Gamma_k) \subset G^2$.

\bigskip
\bigskip
{\bf \S 3. Orthogonal invariants}

\bigskip We keep the notation of the previous sections. In particular, $G$ is a finite group,  $L$ is a $G$-Galois algebra, and $\phi : \Gamma_k \to G$ is
a homomorphism  corresponding to $L$. Let us suppose that the $H^1$-condition is satisfied.

\medskip 
Let $A$ be an {\it orthogonal}  $\sigma^s$-stable central simple factor of $k[G]^s$ (see 1.3), and recall that the center of $A$ is denoted by $E_A$. 
Let us denote by $\langle A \rangle$ the subgroup of  ${\rm Br}(E_A)$ generated by the class of the algebra $A$. Note that since $\sigma_A : A \to A$ is an orthogonal involution, this class has order at most 2, hence $\langle A \rangle$ is a subgroup of  ${\rm Br_2}(E_A)$. 

\medskip
The aim of this section is to define two invariants :
an invariant $c_A(L) \in H^2(k)$ of the $G$--Galois algebra $L$, and an invariant ${\rm clif}_A(q_L) \in {\rm Br_2}(E_A)/ \langle A \rangle$ of the $G$-form $q_L$. We shall
compare these two invariants (cf. Theorem 3.3), and give a necessary condition for the existence of self-dual normal bases (Corollary 3.5). 

\medskip
Let $U^0_A$ be the
connected component of the identity in $U_A$. 
Let $i_A : G \to U_A(E_A)$   be the homomorphism defined in 1.5, and let $\pi : U_A(E_A) \to U_A(E_A)/U^0_A(E_A)$ be the
projection. Since  
$U_A(E_A)/U^0_A(E_A)$  is of order $\le 2$, we have $\pi(i_A(G^2)) = 0$; i.e.
$i_A  (G^2) \subset U_A^0(E_A)$. 

\medskip

Let $\tilde U_A$ be the Spin group of $(A,\sigma)$; note that if ${\rm dim}_k(A) \ge 3$, then $\tilde U_A$ is the
universal cover of $U^0_A$. Let $s : \tilde U_A  \to U^0_A$ be the  covering map.  We have an exact sequence of algebraic groups over $E_A$

$$ 1 \to {\bf Z}/2{\bf Z} \to \tilde U_A  \ {\buildrel s \over \to}  \  U_A^0 \to 1.$$


Let us consider the associated cohomology  exact sequence

$$\tilde U_A(E_A) \ {\buildrel s \over \to} \ U_A^0(E_A) \ {\buildrel \delta \over  \to} \  H^1(E_A).$$

\noindent
{\bf Lemma 3.1.} {\it We have $i_A(G^2) \subset s(\tilde U_A(E_A))$.}

\medskip
\noindent
{\bf Proof.} In view of the above exact sequence, it suffices to prove that $\delta (i_A(G^2))  = 0$. In order to prove this, let us
first assume that $A$ is not split. Then we have $U_A(E_A) = U_A^0(E_A)$ (cf. [K 69], Lemma 1 b, see also [B 94], cor. 2).
Since  $H^1(E_A)$ is a 2-torsion group and since $i_A(G^2) \subset U_A^0(E_A)$, this implies that 
$\delta (i_A (G^2))= 0$, as claimed. Assume now that $A$ is split. Then $U_A$ is the orthogonal group of a quadratic form $q$; let ${\rm sn} : U_A(E_A) \to H^1(E_A)$ be the associated spinor norm, and note that ${\rm sn}$ is a group homomorphism  (see for instance [L 05], Chapter 5, Theorem 1.13). The homomorphism  ${\rm sn}$ depends on the choice of the
quadratic form $q$ with orthogonal group $U_A$, but its restriction to $U^0_A$ does not depend on this choice. Note that
$\delta : U_A^0(E_A) \to H^1(E_A)$ is the restriction of sn to $U_A^0(E_A)$. Therefore for all $g \in G$, we
have $\delta (i_A(g^2)) = {\rm sn} (i_A(g))^2$, and since $H^1(E_A)$ is a 2-torsion group, this implies that  $\delta (i_A(G^2))  = 0$.
This completes the proof of the lemma.

\medskip
Let $H$ be a subgroup of $G^2$. By Lemma 3.1, we have $i_A(H) \subset s(\tilde U_A(E_A))$. Let $$V_A^H = \tilde U_A(E_A) \times_{ U^0_A(E_A) } H$$ be the fibered product
of $s : \tilde U_A(E_A)   \to U^0_A(E_A) $ and $i_A : H  \to U^0_A(E_A) $. 
Therefore we have a central extension

$$1 \to {\bf Z}/2{\bf Z} \to V_A^H  \  {\buildrel p \over \to} \  H \to 1,$$ where 
$p$ is the projection to the factor $H$. Note that the surjectivity of $p$ follows from the fact that by Lemma 3.1 every element of $i_A(H)$ has a preimage in $\tilde U_A(E_A)$.

\medskip
Let us
denote by $$e^H_A \in H^2(H)$$ the cohomology class 
corresponding to the extension $V^H_A$. If $\phi(\Gamma_k) \subset H$, we denote by
$$\phi^* : H^2(H) \to H^2(k)$$ the homomorphism induced by $ \phi : \Gamma_k \to H$.

\medskip
\noindent
{\bf Proposition 3.2.} {\it Let $\psi: \Gamma_k \to G$ be another continuous homomorphism corresponding to the $G$--Galois algebra $L$. Set 
$H_{\phi} = \phi(\Gamma_k)$ and $H_{\psi} = \psi(\Gamma_k)$. Then we have

$$\phi^*(e_A^{H_{\phi}})=  \psi^*(e_A^{H_{\psi}})\ \ {in}  \ \ H^2(k).$$}

\medskip
\noindent
{\bf Proof.} We have $\psi = {\rm Int}(g) \circ \phi$ for some $g \in G$. Note that $i_A(g) \in U_A(E_A)$, and that ${\rm Int}(i_A(g))$ is an automorphism of $U_A^0(E_A)$. 
Any automorphism of $U_A^0(E_A)$ can be lifted to an automorphism of $\tilde U_A(E_A)$; indeed, such a lift exists over a separable closure, and is unique, hence defined
over the ground field.  Let $f : \tilde  U_A(E_A) \to \tilde U_A(E_A)$ be a lift of ${\rm Int}(i_A(g))$. Then $f$ induces an isomorphism $V^{H_{\phi}}_A \to V^{H_{\psi}}_A$, which
sends $H_{\phi}$ to $H_{\psi}$, and is the identity on ${\bf Z}/2{\bf Z} $. This implies that $\phi^*(e_A^{H_{\phi}})=  \psi^*(e_A^{H_{\psi}})\ \ {\rm in}  \ \ H^2(k).$




\medskip
{\bf The invariant $c_A(L)$}

\medskip
We now choose for $H$ the image  $\phi(\Gamma_k)$ of $\Gamma_k$ in $G$, and set $V_A = V_A^H$, $e_A = e_A^H$.
We denote by $c_A(L)$  the class of $\phi^*(e_A^{})$ in $H^2(k)$; Proposition 3.2 shows that this class 
does not depend on the choice of $\phi : \Gamma_k \to G$ defining the $G$--Galois algebra $L$. Since $H^2(k) \simeq {\rm Br}_2(k)$, we can also consider
$c_A(L)$ as an element of ${\rm Br}_2(k)$.

\medskip
Recall that the $G$-trace form $q_L$ determines a 
rank one hermitian form over $(A, \sigma_A)$. We want to relate $c_A(L)$ to the Clifford invariant of this hermitian form. 

\medskip
\medskip
{\bf The invariant ${\rm clif}_A(q_L)$}
\medskip

The map $i_A : H \to U^0_A(E_A)$ 
induces a map of pointed sets  $$i_A : H^1({E_A},H) \to  H^1(E_A,U^0_A).$$
Let  $u^0_A(L)$  be the image of $[\phi_A] \in H^1(E_A,H)$ by
this map. Then the element  $u_A(L)$ defined in 1.5 is the image of $u^0_A(L)$ under the further
composition with the map $H^1(E_A, U^0_A)  \to
H^1(E_A,U_A)$. 


\medskip

Let us consider the exact sequence 
$ 1 \to {\bf Z}/2{\bf Z} \to \tilde U_A  \to U^0_A  \to 1,$ and
let $\delta$ be the
connecting map $H^1(E_A,U^0_A) \to H^2(E_A) \simeq {\rm Br}_2(E_A)$
of the associated cohomology exact sequence. Recall that  $\langle A \rangle$ is the subgroup of  ${\rm Br}_2(E_A)$ generated by the class of the algebra $A$.
The {\it Clifford invariant} of $q_L$  at $A$ is by definition the image  of $\delta(u_A^0(L))$ in $ {\rm Br}_2(E_A)/ \langle A \rangle$. 
Let us denote it by ${\rm clif}_A(q_L)$.

\medskip
\noindent
{\bf Theorem 3.3.} {\it   The image of ${\rm Res}_{E_A/k}(c_A(L))$ in $ {\rm Br}_2(E_A)/ \langle A \rangle$ is equal to  ${\rm clif}_A(q_L)$.}


\bigskip
We need the following lemma :

\medskip
\noindent
{\bf Lemma 3.4.} {\it Let $K$ be a field, let $C$ be a finite group, and let $f : \Gamma_K \to C$ be a continuous homomorphism. Let us denote by $[f] \in H^1(K,C)$
the corresponding cohomology class. Let $$1 \to {\bf Z}/2{\bf Z}  \to V  \to C \to 1$$ be a central extension
with trivial  $\Gamma_K$-action. Let $[e] \in H^2(C)$ be the class of
a 2-cocycle $e : C \times C \to {\bf Z}/2{\bf Z} $ representing this extension. Let 
$\partial : H^1(K,C) \to H^2(K)$ be the connecting map associated to the above
exact sequence, and let $f^* : H^2(C) \to H^2(K)$ be the map induced by $f$. Then 

$$f^*([e]) = \partial ([f]).$$}

\medskip
\noindent
{\bf Proof.} This follows from a direct computation. For all $\sigma, \tau \in \Gamma_K$, we have 
$f^*(e)(\sigma, \tau) = e(f(\sigma),f(\tau)) = x_{f(\sigma)}x_{f(\tau)}x_{f(\sigma \tau)}^{-1}$,
where $x : C \to V$ is a section. On the other hand,  $(\partial f)(\sigma,\tau) =
 x_{f(\sigma)} { }^{f(\sigma)}(x_{f(\tau)}) x^{-1}_{f(\sigma \tau)}$, and this is equal to  
 $ x_{f(\sigma)}x_{f(\tau)}x_{f(\sigma \tau)}^{-1},$ since the action of $\Gamma_k$
 on $V$ is trivial. 

\medskip
\noindent
{\bf Proof of Theorem 3.3.} 
Let $\partial : H^1(E_A,H) \to H^2(E_A)$ be the connecting map of the cohomology exact sequence associated to 
the exact sequence 
$$1 \to {\bf Z}/2{\bf Z} \to V_A  \to H \to 1$$
with all the groups having trivial $\Gamma_{E_A}$-action. 
Recall that  $\phi_A  : \Gamma_{E_A} \to \Gamma_k \to H$ is the composition of 
$\phi : \Gamma_k \to H$ with the inclusion of $\Gamma_{E_A}$ into $\Gamma_k$. By  Lemma 3.4 we have
$\partial ([\phi_A] ) = \phi_A^*(e_A) =
{\rm Res}_{E_A/k} ( \phi^*(e_A)) = {\rm Res}_{E_A/k} (c_A(L))$. In view of the commutative
diagram of $\Gamma_{E_A}$-groups

$$\matrix {1 & \to & {\bf Z}/2{\bf Z} & \to & \tilde U_A (k_s) & \to & U_A^0(k_s)& \to & 1 \cr  {} & {} & \uparrow  & &  \uparrow & & \uparrow \cr 1 & \to & {\bf Z}/2{\bf Z}  & \to & V_A &  \to & H  & \to & 1  \cr}$$
we have $\delta(u_A^0(L)) = \partial ([\phi_A])$. Therefore we obtain 
$ {\rm Res}_{E_A/k} (c_A(L)) =  \delta(u_A^0(L))$. Since  the class of  $ \delta(u_A^0(L))$ in  ${\rm Br}_2(E_A)/\langle A \rangle $ is equal to ${\rm clif}_A(q_L)$ by definition, this completes the proof of the theorem.

\medskip
\noindent
{\bf Proposition  3.5.} {\it If $L$ has a self-dual normal basis over $k$, then 
${\rm Res}_{E_A/k}(c_A(L))$
is trivial in ${\rm Br}_2(E_A)/\langle A \rangle $.}

\bigskip
\noindent
{\bf Proof.} Since $L$ has a self-dual normal basis over $k$, the class $u_A(L)$ corresponds 
to the class of the rank one unit hermitian form 
$\langle 1 \rangle$ in $H^1(E_A,U_A)$. As  $\langle 1 \rangle$ corresponds to the trivial cocycle in
$Z^1(E_A,U^0_A)$,  its Clifford invariant is trivial, in other words,  ${\rm clif}_A(q_L)$
is trivial. By Theorem 3.3 
the 
image of ${\rm Res}_{E_A/k}(c_A(L))$ in $ {\rm Br}_2(E_A)/ \langle A \rangle$ is equal to  ${\rm clif}_A(q_L)$, hence the proposition is proved. 


\bigskip
We conclude this section with an example where $c_A(L) \not = 0$, but ${\rm Res}_{E_A/k}(c_A(L)) = 0$ (and hence ${\rm clif}_A(q_L) = 0$) :

\bigskip
\noindent
{\bf Example 3.6.} Let $G = A_5$, the alternating group, and assume that $k = {\bf Q}$. Let $A$ be a factor of $k[G]$ corresponding to a degree 3 orthogonal
representation of $G$; then $A = M_3(E_A)$ with $E_A = {k}(\sqrt {5})$, and the involution $\sigma_A$ is induced by the unit form $\langle 1,1,1 \rangle$.
Let $\epsilon \in G$ be  a product of two disjoint transpositions. 

\medskip
Let $z \in k^{\times}$, and let $\psi : \Gamma_k \to \{1,\epsilon \}$ be the corresponding quadratic character. Let $\phi : \Gamma_k \to G$ be given by $\phi = \iota \circ \psi$, where $\iota : \{1,\epsilon \} \to G$ is the inclusion.
Let $L$ be the $G$-Galois algebra corresponding to $\phi$.  Set $H = \{1,\epsilon \}$, and note that the image of $\phi$ is contained in $H$. 
Set $N = k[X]/(X^2-z)$; then we have $L = {\rm Ind}^G_H(N)$.

\medskip
Note that  $\epsilon$ lifts to an element of order 4 in $\tilde A_5$,
hence also in $\tilde U_A(E_A)$. Therefore the extension $1 \to {\bf Z}/2{\bf Z} \to V_A^H  \to  H \to 1$ is not trivial; the group $V_A^H$ is
cyclic of order 4.  Recall that $e_A$ is the class of this extension in $H^2(H)$; hence $e_A$ is the only non-trivial element of $H^2(H)$. By definition,
we have $c_A(L) = \phi^*(e_A)$, and this is equal to the cup product $(z)(z) = (-1)(z)$ in $H^2(k)$. 

\medskip
Set $z = 5$. Then $c_A(L) = (-1)(5)$ is not trivial in $H^2(k)$. On the other hand, since $E_A = {k}(\sqrt {5})$, we have ${\rm Res}_{E_A/k}(c_A(L)) = 0$ in $H^2(E_A)$. Note
that the subgroup $\langle A \rangle$ of ${\rm Br}_2(E_A)$ is trivial, and
recall
that ${\rm clif}_A(q_L) = {\rm Res}_{E_A/k}(c_A(L))$ in ${\rm Br}_2(E_A) \simeq H^2(E_A)$ by Theorem 3.3; therefore we have ${\rm clif}_A(q_L) = 0$.

\bigskip
\bigskip
{\bf \S 4. Unitary invariants}

\bigskip We keep the notation of the previous sections : $G$ is a finite group, $L$ is a $G$--Galois algebra, and $\phi : \Gamma_k \to G$ is
a homomorphism associated to $L$. We
suppose that the $H^1$-condition is satisfied by $\phi : \Gamma_k \to G$, hence $\phi (\Gamma_k)$ is a subgroup of $G^2$. 
Let $A$ be a {\it unitary}  $\sigma^s$-stable central simple factor of $k[G]^s$ (see 1.3). 
We denote by $F_A$ be the center of $A$; note that $F_A$
is a quadratic extension of $E_A$. 

\medskip
Using the same strategy as in \S 3, we first define an element of $H^2(k)$ which is an invariant of the $G$-Galois algebra $L$. We then consider the hermitian form $h_A$ over
$(A,\sigma)$ determined by $q_L$, and recall the definition of the discriminant of this form, thereby obtaining an element of ${\rm Br}_2(E_A)$. This 
is an invariant of the hermitian form $h_A$, and hence also of the  $G$--form $q_L$. We then show that the restriction of the first invariant to $H^2(E_A)$ is equal to the second one (see Theorem 4.5).

\medskip
We start by recording some facts from Galois cohomology.


\medskip
Let $E$ be a field of characteristic $\not =2$, and let $E_s$ be a separable closure of $E$. 
Let $F$ be a quadratic extension of $E$, let  $x \mapsto \overline x$
the non-trivial automorphism of $F$ over $E$, and let $F^{{\times}1}$ be the subgroup
of $F^{\times}$ consisting of the $x \in F$ such that $x \overline x = 1$. 
Let ${\rm N }: F \to E$, given by ${\rm N}(x) = x \overline x$, be the norm map.
We denote by $[F]$ the
 class of the quadratic extension $F/E$ in $H^1(E)$. For all $x \in E^{\times}$, we denote by $(x)$ the class of $x$ in $E^{\times}/E^{{\times}2}$, and by
 $[x]$ the class of $x$ in $E^{\times}/{\rm N}(F^{\times})$.

\medskip

\medskip
\noindent
{\bf Lemma 4.1.} {\it {\rm (a)} The connecting homomorphism $E^{{\times}}  \to H^1(E,{\rm R}^1_{F/E} {\bf G}_m)$
associated to the exact sequence \ $1 \to {\rm R}^1_{F/E} {\bf G}_m \to {\rm R}_{F/E} {\bf G}_m   \ { \buildrel {\rm N} \over \to}  \ {\bf G}_m   \to 1$ \  induces an isomorphism
$\alpha : E^{{\times}}/{\rm N}(F^{\times})  \to H^1(E,{\rm R}^1_{F/E} {\bf G}_m)$.

{\rm (b)}
Let $x \in E^{\times}$, and let $f_x : \Gamma_E \to {\rm R}^1_{F/E} {\bf G}_m (E_s)$ 
be defined by $f_x(\gamma) = y^{-1} \gamma (y)$,
where $y \in (F  \otimes_E  E_s)^{\times}$  is such that ${\rm N}(y) = x$. Then we have $\alpha((x)) = [f_x]$.}

\medskip
\noindent
{\bf Proof.} (a) follows from Hilbert's theorem 90, and (b) from the definition of the connecting homomorphism. 

\medskip
From now on, we identify $E^{{\times}}/{\rm N}(F^{\times})$ and $H^1(E,{\rm R}^1_{F/E} {\bf G}_m )$ via the isomorphism $\alpha$.



\medskip
 \noindent
 {\bf Lemma 4.2.} {\it Let $1 \to  {\bf Z}/2{\bf Z} \to {\rm {\rm R}}^1_{F/E} {\bf G}_m \  {\buildrel s \over \to} \
{\rm R}^1_{F/E} {\bf G}_m \to 1$ be the exact sequence of linear algebraic groups with $s$ the squaring map.
 Let $\delta : H^1(E,  {\rm R}^1_{F/E} {\bf G}_m) \to H^2(E)$  be the connecting homomorphism associated to this exact sequence. Identifying $H^1(E, {\rm R}^1_{F/E} {\bf G}_m)$
 with $E^{\times}/{\rm N}(F^{\times})$ via $\alpha$, 
 we have
 $$\delta ([x]) = (x)[F] \in H^2(E)$$
 for all $x \in E^{\times}$, where $(x)[F]$ denotes the cup product of $(x),[F] \in H^1(E)$.}
 
 \medskip
 \noindent  
 {\bf Proof.} 
 A 2-cocycle 
 associated to $(x)[F] \in H^2(E)$ is given by $f(\sigma,\tau)$ such that
 $f(\sigma,\tau) = 1$ if the restriction of $\sigma$ to $E(\sqrt x)$ is the identity, or if
 the restriction of $\tau$ to $F$ is the identity, and $f(\sigma,\tau) = - 1$ otherwise.
 Let us check that the cohomology class of $f$ in $H^2(E)$ is equal to $\delta([x])$.
 Let $y \in (F  \otimes_E  E_s)^{\times}$ be such that ${\rm N}_{F  \otimes_E  E_s/E_s}(y)  = y \overline y = x$. A
 1-cocycle $g$ in $Z^1(E, {\rm R}^1_{F/E} {\bf G}_m)$ associated to $[x]$ is given by
 $g(\sigma) = y^{-1} \sigma (y)$ for $\sigma \in \Gamma_E$. For all $\tau \in \Gamma_E$, set $z_{\tau} =
 y^{-1} \sqrt x$ if the restriction of $\tau$ to $F$ is not the identity, and 
 $z_{\tau} = 1$ otherwise. Then ${\rm N}_{F  \otimes_E  E_s/E_s}(z_{\tau})  = z_{\tau} \overline {z_{\tau}} = (y^{-1} \sqrt x) (\overline y ^{-1} \sqrt x)$
 if the restriction of $\tau$ to $F$ is not the identity. Since $y \overline y = x$, we have
 $z_{\tau} \in {\rm R}^1_{F/E} {\bf G}_m (E_s)$. Further,  $s(z_{\tau}) = y^{-2} x = y^{-1} \tau(y)$
 if the restriction of $\tau$ to $F$ is not the identity, and $s(z_{\tau}) = 1 = y^{-1} \tau (y)$
 otherwise. Thus $\delta(g)(\sigma,\tau) = z_{\sigma}{}^{\sigma} (z_{\tau}) z_{\sigma \tau}^{-1}$.
 It is straightforward to check that $\delta (g) (\sigma,\tau) = 1$ if the restriction
 of $\sigma$ to $E(\sqrt x)$ is the identity, or the restriction of $\tau$ to $F$ is
 the identity, and that  $\delta (g) (\sigma,\tau) = -1$ otherwise. This is precisely the cocycle $f$, hence we have $\delta([x]) = (x)[F]$
 in  
 $H^2(E)$. This concludes the proof of the lemma.
 
 
\medskip
\noindent
{\bf Lemma 4.3.}
{\it  We have an injective homomorphism $E^{{\times}}/{\rm N}(F^{\times}) \to {\rm Br_2}(E)$ defined by $[x] \mapsto (x,F/E)$.}

\medskip
\noindent
{\bf Proof.} Indeed, the class of the quaternion algebra $(x,F/E)$ is trivial in   ${\rm Br_2}(E)$  if and only if $x \in {\rm N}(F^{\times})$.

\medskip

We now define an invariant $d_A(L) \in H^2(k,  {\bf Z}/2{\bf Z})$ of the $G$-Galois algebra $L$. 

\medskip

{\bf The invariant $d_A(L)$}

\medskip
Recall that $F_A^{{\times}1}$ is  the subgroup
of $F_A^{\times}$ consisting of the $x \in F_A$ such that $x \sigma_A (x) = 1$; 
in other words, $F_A^{{\times}1} = {\rm R}^1_{F_A/E_A}{\bf G}_m(E_A)$. We denote by $s : {\rm R}^1_{F_A/E_A}{\bf G}_m \to {\rm R}^1_{F_A/E_A}{\bf G}_m$ the squaring map,
and by $n : U_A \to {\rm R}^1_{F_A/E_A}{\bf G}_m$ the reduced norm. 
Recall that $i_A : G \to U_A(E_A)$ is
the homomorphism defined in 1.5; we have $n(i_A(G^2)) \subset s(F_A^{{\times}1}).$

\medskip
Let $H$ be a subgroup of $G^2$.
Let $V^H_{A} = F^{{\times}1}_{A} \times_{F^{{\times}1}_{A}} H$ be the fibered product
of $s : F^{{\times}1}_{A} \to F^{{\times}1}_{A}$ and $n \circ i_A : H \to F^{{\times}1}_{A}$. Then
the sequence 

$$1 \to {\bf Z}/2{\bf Z} \to V_A^H   \to H \to 1$$ is exact. Note that the surjectivity follows from the fact that $n(i_A(H))  \subset s(F_A^{{\times}1}).$
Therefore $V^H_A$ is a central extension
of $H$ by ${\bf Z}/2{\bf Z}$. Recall that the $H^1$-condition implies that $\phi (\Gamma_k) \subset G^2$. 

\medskip
\noindent
{\bf Proposition 4.4.} {\it Let $\psi: \Gamma_k \to G$ be another continuous homomorphism corresponding to the $G$--Galois algebra $L$. Set 
$H_{\phi} = \phi(\Gamma_k)$ and $H_{\psi} = \psi(\Gamma_k)$. Then we have

$$\phi^*(e_A^{H_{\phi}})=  \psi^*(e_A^{H_{\psi}})\ \ { in}  \ \ H^2(k).$$}

\medskip
\noindent
{\bf Proof.} We have $\psi = {\rm Int}(g) \circ \phi$ for some $g \in G$. 
The map $F_A^{\times 1} \times_{F^{\times 1}} H_{\phi} \to F_A^{\times 1} \times_{F^{\times 1}} H_{\psi}$, 
given by $(x,y) \to (x,gyg^{-1})$, gives rise to
an isomorphism $V_A^{H_{\phi}} \to V_A^{H_{\psi}}$ that is the identity on ${\bf Z}/2{\bf Z}$ and sends $H_{\phi}$ to $H_{\psi}$.
This implies that $\phi^*(e_A^{H_{\phi}})=  \psi^*(e_A^{H_{\psi}})\ \ {\rm in}  \ \ H^2(k).$

\medskip 
We now choose for $H$ the image  $\phi(\Gamma_k)$ of $\Gamma_k$ in $G$, and set $V_A = V_A^H$, $e_A = e_A^H$.

\medskip
\noindent
{\bf Notation.}
Let us denote by $d_A(L)$  the class of $\phi^*(e_A)$ in $H^2(k)$; Proposition 4.4 shows that  this class is independent of the choice of $\phi : \Gamma_k \to G$ defining the $G$--Galois algebra $L$.

\bigskip

We define the discriminant of the $G$-form $q_L$ at $A$, and compare it with the cohomology class $d_A(L)$.

\medskip

{\bf The invariant  ${\rm disc}_A(q_L)$}



\medskip
Recall that composing  $\phi_A : \Gamma_{E_A} \to H$ with the
map $i_A : H \to U_A(k_s)$ we obtain a $1$-cocycle $\Gamma_{E_A} \to U_A(k_s)$, the class of which in $H^1(E_A,U_A)$
is $u_A(L)$. The reduced norm $n : U_A \to {\rm R}^1_{F_A/E_A}{\bf G}_m$  induces
a map $n : H^1(E_A,U_A) \  {\buildrel {} \over \to} \ E_A^{{\times}}/{\rm N}(F_A^{\times}).$



\medskip
\noindent
{\bf Notation.}  Set ${\rm disc}_A(q_L) = (n(u_A(L)),F_A/E_A)$ in ${\rm Br}_2(E_A)$. 


\medskip Note that this is well--defined by Lemma 4.3. Since we have ${\rm Br_2}(E) \simeq H^2(E_A)$, we can also consider ${\rm disc}_A(q_L)$ as an element of  $H^2(E_A) $.
Then ${\rm disc}_A(q_L)$ is given by the cup product $ n(u_A(L)).[F_A]$ in $H^2(E_A) $.
This invariant is related to the previously defined invariant $d_A(L)$ as follows : 

\medskip
 \noindent
 {\bf Theorem 4.5.} {\it We have
 ${\rm disc}_A(q_L) = {\rm Res}_{E_A/k}(d_A(L))$  in $H^2(E_A)$.}

 \medskip
 \noindent
 {\bf Proof.} Let $\partial : H^1(E_A,H) \to H^2(E_A)$ be the connecting map of
the exact sequence 
$$1 \to {\bf Z}/2{\bf Z} \to V_A  \to H \to 1$$
with all the groups having trivial $\Gamma_{E_A}$-action. 
 By  Lemma 3.4 we have
$$\partial ([\phi_A] ) = \phi_A^*(e_A) =   {\rm Res}_{E_A/k} ( \phi^*(e_A)) = {\rm Res}_{E_A/k} (d_A(L)).$$
 We have the commutative diagram
 $$\matrix {1 & \to &  {\bf Z}/2{\bf Z} & \to & V_A & \to & H & \to & 1 \cr  {} & {} & \downarrow & &  \downarrow & & \downarrow \cr 1 & \to &   {\bf Z}/2{\bf Z}   & \to & 
 {\rm R}^1_{F_A/E_A}  {\bf G}_m (k_s) &  {\buildrel s \over \to} & {\rm R}^1_{F_A/E_A} {\bf G}_m (k_s)    & \to & 1  \cr}$$
 where the second vertical map is the projection on the first factor, 
and the third one is 
 $H \  {\buildrel {i_A} \over \to} \ U_A(E_A) \  {\buildrel n \over \to} \  {\rm R}^1_{F_A/E_A} {\bf G}_m (E_A)$.
 
 \medskip
 Let $\delta : H^1(E,  {\rm R}^1_{F_A/E_A} {\bf G}_m) \to H^2(E_A)$  be the connecting homomorphism associated to the exact sequence
 
 $$1  \to    {\bf Z}/2{\bf Z}    \to 
 {\rm R}^1_{F_A/E_A}  {\bf G}_m  \  {\buildrel s \over \to}  \  {\rm R}^1_{F_A/E_A} {\bf G}_m  \to  1.$$ By the commutativity of the above diagram, we
 have $\delta ([n(u_A(L))]) = \partial ([\phi_A] )$. 
 Hence we have ${\rm Res}_{E_A/k}(d_A(L)) =  \delta ([n(u_A(L)]).$ We have $ \delta ([n(u_A(L))]) =  (n(u_A(L))).[F_A]$ by Lemma 4.2
 and hence ${\rm Res}_{E_A/k}(d_A(L)) = {\rm disc}_A(q_L)$, as claimed.
 
 \medskip
\noindent
{\bf Lemma 4.6.}
{\it If $q_L$ corresponds to the hermitian form $\langle z_A \rangle$ over $(A,\sigma_A)$, then we have $${\rm disc}_A(q_L) = (n(z_A),F_A/E_A) \ \  {in}  \ \ {\rm Br}_2(E_A).
$$}

\noindent
{\bf Proof.} Set $z = z_A$. Let $z = w \sigma_A(w)$ with $w \in A \otimes_{E_A} k_s$. The cocycle $\tau \mapsto w^{-1} \tau (w)$ represents the class of the hermitian form $\langle z \rangle$ in $H^1(E_A,U_A)$. Let us denote this class by 
$u_z \in H^1(E_A,U_A)$, and note that we have $u_z = u_A(L)$ by definition. 
The cocycle $\tau \mapsto n(w)^{-1} \tau (n(w))$ represents  the class $n(u_z) \in H^1(E_A,{\rm R}^1_{F_A/E_A} {\bf G}_m)$.
By Lemma 4.1 this class is
mapped by $\alpha^{-1}$ to $[n(z)]  \in E_A^{{\times}}/{\rm N}(F_A^{{\times}})$. Therefore we have $(n(z),F_A/E_A) = (n(u_A(L)),F_A/E_A)  = {\rm disc}_A(q_L)$, as claimed.


\medskip
\noindent
{\bf  Proposition 4.7.} {\it If  $L$ has a self-dual normal basis over $k$, then
${\rm Res}_{E_A/k}(d_A(L))$
is trivial in ${\rm Br}_2(E_A)$.}

\medskip
\noindent
{\bf Proof.} Since $L$ has a self-dual normal basis,  $q_L$ corresponds to the hermitian form $\langle 1 \rangle$ over $(A,\sigma_A)$.
By Lemma 4.6 this implies that  ${\rm disc}_A(q_L) $ is trivial.
Since by Theorem 4.5 we have ${\rm disc}_A(q_L) = {\rm Res}_{E_A/k}(d_A(L))$, the Proposition
is proved.

\medskip
\noindent
{\bf Remark.} There are examples where $d_A(L) \not = 0$ but ${\rm Res}_{E_A/k}(d_A(L))= 0$ (hence  also ${\rm disc}_A(q_L) = 0$); see for instance
Example 5.2. {\rm (i)}.

 \bigskip
 \bigskip
 {\bf \S 5. Self-dual normal bases}
 
 \bigskip
 We keep the notation of the previous sections. In particular, $G$ is a finite group, $L$ is a $G$-Galois algebra over $k$, and $\phi : \Gamma_k \to G$ is
 a homomorphism associated to $L$.   We now apply the results of the previous sections to give necessary conditions for the existence of a
self-dual normal basis, and to show that these are also sufficient when $k$ has cohomological dimension $\le 2$,
see Proposition 5.1 and Theorem 5.3. 

\medskip
Putting together the results of \S 2 - \S 4, we have the following :
 
 \medskip
 \noindent
{\bf Proposition 5.1.} {\it Suppose that $L$ has a self-dual normal basis over $k$. Then the $H^1$-condition is satisfied, and



\medskip

{\rm (i)} For all orthogonal $\sigma^s$-stable central simple factors $A$ of $k[G]^s$, we have

\medskip
\centerline {${\rm Res}_{E_A/k}(c_A(L)) = 0$
 \ in ${\rm Br}_2(E_A)/\langle A \rangle$.}

\medskip
{\rm (ii)} For all unitary $\sigma^s$-stable central simple factors $A$ of $k[G]^s$, we have
$${\rm Res}_{E_A/k}(d_A(L)) = 0 \  \ { in} \  {\rm Br}_2(E_A).$$}

\medskip
\noindent
{\bf Proof.} This follows from Propositions 2.1, 3.5 and 4.7.

\medskip
\noindent
{\bf Example 5.2.} {\rm (i)} The aim of this example is to reinterpret and complete Exemple 10.2. of [BSe 94] using the results of the present paper.
Assume that $G$ is cyclic of order 8, and let $s$ be a generator of $G$; let $\epsilon = s^4$ be the element of order 2 of $G$. Let $z \in k^{\times}$, and let $\sigma : \Gamma_k \to \{1,\epsilon \}$ be the corresponding quadratic character. Let $\phi : \Gamma_k \to G$ be given by $\phi = \iota \circ \sigma$, where $\iota : \{1,\epsilon \} \to G$ is the inclusion.
Let $L$ be the $G$-Galois algebra corresponding to $\phi$.  Set $H = \{1,\epsilon \}$, and note that the image of $\phi$ is contained in $H$. 
Set $N = k[X]/(X^2-z)$; then we have $L = {\rm Ind}^G_H(N)$.  Set $A = k[X]/(X^4 + 1)$, and let us write $k[G] = A' \times A$. It is easy to see that the image of $H$ in $A'$ is trivial.
The involution $\sigma_A$ sends the class of $X$ to the class of $X^{-1}$. If $k$ contains the 4th roots of unity, then $A$ is
a product of two factors exchanged by the involution, hence there $k[G]$ has no involution invariant factor in which the image of $H$ is non trivial. In this case, $L$ has a self-dual normal basis. Assume that $k$ does not contain the 4th roots of unity. Then $A$ is a field; we have $F_A = A$, and $E_A = k[X]/(X^2 - 2)$. Note that $A$ is unitary. 
We have $i_A(\epsilon) = -1$, hence $i_A(H) = \{1,-1\}$.

\medskip
Let $i \in F_A$ be a primitive 4th root of unity. By the definition of the extension $$1 \to {\bf Z}/2{\bf Z} \to V_A  \to H \to 1$$ (cf. \S 4), we see that
$V_A = \{(1,1), (-1,1), (i,\epsilon), (-i, \epsilon) \}$, a cyclic group of order 4. Recall that $e_A$ is the class of this extension in $H^2(H)$;
hence $e_A$ is the only non-trivial element of $H^2(H)$. We have $d_A(L) = \phi^*(e_A) = (z,z) = (z,-1)$, and 
${\rm Res}_{E_A/k}(d_A(L)) = (z,-1)_E = (z,F_A/E_A)$. Therefore we have 

\medskip
\centerline {$d_A(L)  = 0 \iff$ $z$ is a sum of two squares in $k$,} and

\centerline {${\rm Res}_{E_A/k}(d_A(L))  = 0 \iff$ $z$ is a sum of two squares in $E_A = k(\sqrt 2)$.}

\medskip It is easy to find examples where $d_A(L)  \not = 0$ and ${\rm Res}_{E_A/k}(d_A(L))  = 0$; for instance, we can take $k = {\bf Q}$ and $z = 3$.

\medskip By Proposition 5.1 the existence of a self-dual normal basis implies that we have ${\rm Res}_{E_A/k}(d_A(L))  = 0$. On the other hand, in [BSe 94], Exemple 10.2 it is checked
by direct computation that if $z$ is a sum of two squares in $k(\sqrt 2)$, then $L$ has a self-dual normal basis. Hence we have

\medskip

\centerline {$L$ has a self-dual normal basis over $k  \iff$ {\it z} is a sum of two squares in $k(\sqrt 2)$.}

\medskip
Moreover,  $z$ is a sum of two squares in $k(\sqrt 2)$ if and only if $z$ is a sum of 4 squares in $k$. Indeed, $z$ is a sum of two squares in $k(\sqrt 2)$ $\iff$ the quadratic form
$\langle 1,1, -z,-z \rangle$ represents 0 over $k(\sqrt 2)$  $\iff$ $\langle 1,1, -z,-z \rangle$ represents -2 over $k$ $\iff$ the quadratic form $\langle 1,1, -z,-z \rangle \otimes \langle 1,2 \rangle$ represents zero over $k$ $\iff$ $\langle 1,1, -z,-z \rangle \otimes \langle 1,1 \rangle$ represents zero over $k$ $\iff$  $\langle 1,1,1,1 \rangle$ represents $z$ over $k$. Hence we get

\medskip

\centerline {$L$ has a self-dual normal basis over $k  \iff$ {\it z} is a sum of 4 squares in $k$.}

\medskip
{\rm (ii)} Assume that $G = D_4$, the dihedral group of order 8.  Then a $G$-Galois algebra $L$ has a self-dual normal basis if and only if 
either $L$ is split or $L = {\rm Ind}^G_H(N)$ with $H$ of order 2, and $N = k[X]/(X^2-z)$ for some $z \in k^{\times}$ such that z is a sum of two squares in $k$.

\medskip
Indeed, let $\phi : \Gamma_k \to G$ be a homomorphism associated to $L$. Note that
$G^2$ is of order 2, hence the $H^1$-condition holds if and only if the image of $\phi$ is of order 1 or 2; in other words, $L$ is split, or induced from a group of order $2$. 
Let  $H = \{1,\epsilon \}$, and assume that $L = {\rm Ind}^G_H(N)$, with $N = k[X]/(X^2-z)$ for some $z \in k^{\times}$. 
\medskip

The group $G$ has one degree 2 and four degree 1 orthogonal representations. Since the $H^1$-condition holds, the image of $G$ is trivial in the factors of $k[G]$ corresponding to the degree 1 representations. Let $A = M_2(k)$, and let $\sigma_A$ be the involution induced by the 2-dimensional unit form; then the factor of $k[G]$ corresponding to the degree 2 orthogonal representation of $G$
is equal to $A$. 

\medskip
If $k$ contains the $4$-th roots of unity, then $U_A^0 =  \tilde U_A =  {\bf G}_m$. If $k$ does not contain the $4$-th roots of unity, then
$U_A^0 =  \tilde U_A = R^1_{K/k}{\bf G}_m$, where $K = k[X]/(X^2 + 1)$. In both cases, $s : \tilde U_A \to U_A^0$ is the squaring map. Using this, we see that
the extension $1 \to {\bf Z}/2{\bf Z} \to V_A  \to H \to 1$ is non-trivial, and that 
$c_A(L) = (z,-1)$. Therefore 

\medskip

\centerline {$L$ has a self-dual normal basis over $k  \iff$ {\it z} is a sum of two squares in $k$.}

\medskip
{\rm (iii)} Let $G = A_4$, the alternating group of order 12, and assume for simplicity that $char (k) \not = 3$ and that $k$ contains the third roots of unity. Then $k[G] = k \times k \times k \times M_3(k)$, where the first factor corresponds to the unit representation, the second and the third to the two representations of degree 1 with image of order 3, and the fourth
one to the irreducible representation of degree 3. Let $A = M_3(k)$ be the fourth factor, and note that  the restriction of $\sigma$ to $A$ is induced by the 3-dimensional unit form. The extension $1 \to {\bf Z}/2{\bf Z} \to V_A  \to G \to 1$ defined in \S 3 is

$$1 \to {\bf Z}/2{\bf Z} \to \tilde A_4  \to A_4 \to 1,$$ corresponding to the unique non-trivial element $e \in H^2(A_4)$  (see [Se 84], 2.3). 
Let $L$ be a $G$-Galois algebra, and note that the $H^1$-condition
is satisfied, since $G$ has no quotient of order 2. Let $E$ be the subalgebra of $L$ fixed by the subgroup $A_3$ of $G = A_4$; then $E$ is an \'etale algebra of rank 4. Let $\phi : \Gamma_k \to A_4$ be
a homomorphism corresponding to $L$.
By [Se 84], Theorem 1 we have $\phi^*(e) = w_2(q_E)$,
the Hasse-Witt invariant of the quadratic form $q_E$; hence
the
invariant $c_A(L)$ is equal to $w_2(q_E)$. 
Let $q_A(L)$ be the 3-dimensional quadratic form
corresponding to the cohomology class $u_A(L)$. Then $q_E \simeq q_A(L) \oplus \langle 1 \rangle$, and it is easy to check that
$q_A(L) \simeq \langle 1,1,1\rangle$
$\iff$ $w_2(q_E) = 0$, hence $u_A(L) = 0$ $\iff$ $w_2(q_E) = 0$. Therefore we have 


\medskip
\centerline {$L$ has a self-dual normal basis over $k  \iff$  $w_2(q_E) = 0$,} 

\medskip 
\noindent recovering a result of [BSe 94] (see [BSe 94], Exemple 1.6). 




\bigskip The case of cyclic groups of order a power of 2 is further developed in \S 6; we now look at fields of low cohomological dimension. 
Recall that  the {\it 2-cohomological dimension} of $\Gamma_k$, 
denoted by ${\rm cd}_2(\Gamma_k)$, is the smallest integer $d$ such that $H^i(k,C) = 0$ for all $i > d$ and for every
finite 2-primary $\Gamma_k$-module $C$. For fields of cohomological dimension $\le 1$, the question of existence of self-dual normal bases is settled in [BSe 94], 2.2. 

\medskip
\noindent
{\bf Theorem 5.3.} {\it Assume that 
${\rm cd}_2(\Gamma_k) \le 2$. 
Then $L$ has a self-dual normal basis over $k$ if and only if the $H^1$-condition is satisfied, and 
the conditions {\rm (i)} and  {\rm (ii)}  
below hold  :

\medskip

{\rm (i)} For all orthogonal $\sigma^s$-stable central simple factors $A$ of $k[G]^s$, we have

\medskip
\centerline {${\rm Res}_{E_A/k}(c_A(L)) = 0$ \
in  ${\rm Br}_2(E_A)/\langle A \rangle$.}

\medskip
{\rm (ii)} For all unitary $\sigma^s$-stable central simple factors $A$ of $k[G]^s$, we have
$${\rm Res}_{E_A/k}(d_A(L)) = 0 \  \  { in} \   {\rm Br}_2(E_A).$$}
\noindent
{\bf Proof.} If $L$ has a self-dual normal basis over $k$, then by Proposition 5.1 the $H^1$-condition, as well as conditions {\rm (i)} and {\rm (ii)}  
are satisfied. Conversely, let us assume that the $H^1$-condition, as well as
conditions  {\rm (i)} and {\rm (ii)} 
 hold. Since the $H^1$-condition holds, we can define $c_A(L)$ and $d_A(L)$, cf. \S 3 and \S 4. By Theorems 3.3
and 4.5  we have ${\rm clif}_A(q_L) = {\rm Res}_{E_A/k}(c_A(L))$  and  ${\rm disc}_A(q_L) = {\rm Res}_{E_A/k}(d_A(L))$.  Therefore, conditions {\rm (i)}  and {\rm (ii)} 
imply that ${\rm clif}_A(q_L)$ is trivial for all orthogonal factors $A$, and ${\rm disc}_A(q_L)$ is trivial for all unitary factors $A$.  
Let us prove that $L$ has a self-dual normal basis over $k$. Let us denote by  $h_A$ the hermitian form over $(A,\sigma_A)$ corresponding to $u_A(L)$. It is enough to show that for all factors $A$,  the class $u_A(L)$ is trivial; this is equivalent
with saying that the hermitian form $h_A$ is isomorphic to the unit form $1_A$ over $(A,\sigma_A)$. By Witt
cancellation (see for instance [BPS 13],  Theorem 2.5.2) this in turn is equivalent to saying that $h_A \oplus - 1_A$ is hyperbolic. Let us prove this successively for symplectic,
orthogonal and unitary characters. 

\medskip
Assume first that $A$ is symplectic. Then by [BP 95], Theorem 4.3.1 every even dimensional non-degenerate hermitian form over a central simple algebra with involution
is hyperbolic. This implies that $h_A \oplus - 1_A$ is hyperbolic. Assume now that $A$ is orthogonal, and note that the $H^1$-condition implies that $u_A(L)$ is
the image of a class $u^0_A(L)$ of $H^1(E_A,U_A^0)$. This implies that $h_A$ has trivial discriminant. As we saw above,  ${\rm clif}_A(q_L)$ is trivial, hence
the form $h_A \oplus - 1_A$ has trivial Clifford invariant. By [BP 95], Theorem 4.4.1. every even dimensional  non-degenerate hermitian form over a central
simple algebra having trivial discriminant and trivial Clifford invariant is hyperbolic, hence $h_A \oplus - 1_A$ is hyperbolic. Assume finally that $A$ is a unitary
character. We have seen above that ${\rm disc}_A(q_L)$ is trivial, therefore the form  $h_A \oplus - 1_A$ has trivial discriminant. By [BP 95], Theorem 4.2.1. 
every even dimensional  non-degenerate hermitian form over a central
simple algebra having trivial discriminant is hyperbolic, hence $h_A \oplus - 1_A$ is hyperbolic.

\medskip
This implies that $L$ has a self-dual normal basis over $k$, hence the theorem is proved. 

\medskip
Recall that $\phi : \Gamma_k \to G$ is a homomorphism associated to the $G$-Galois algebra $L$, and that for all  $x \in H^n (G)$, we denote by $x_L$ the image of $x$ by
$\phi^* : H^n(G) \to H^n(k)$. Let $H = \phi(\Gamma_k)$. For $n =2$, we also need the image of $x$ by the homomorphism $\phi^* : H^n(H) \to H^n(k)$;
we denote this image by $x^H_L$.

\medskip
\noindent
{\bf Corollary 5.4.} {\it Assume that ${\rm cd}_2(\Gamma_k) \le 2$, that the $H^1$-condition is satisfied, 
and that  we have  $x^H_L = 0$ for all $x \in H^2(H)$.
Then $L$ has a self-dual normal basis over $k$.}

\medskip
\noindent
{\bf Proof.} This follows immediately from Theorem 5.3. Indeed, the  $H^1$-condition is satisfied by hypothesis. Moreover,
the classes $c_A(L)$ and $d_A(L)$ are by definition in the image of $\phi^* : H^2(H) \to H^2(k)$,
hence the hypothesis $x^H_L = 0$ for all $x \in H^2 (H)$ implies that $c_A(L) = 0$ 
for all orthogonal factors $A$, and $d_A(L) = 0$ for all unitary factors $A$. Therefore conditions  {\rm (i)} and {\rm (ii)} of Theorem 5.3 are
satisfied, and hence $L$ has a self-dual normal basis over $k$. 

\medskip
\noindent
{\bf Remarks.}
{\rm (i)} Corollary 5.4 suggests the following question : Assume that ${\rm cd}_2(\Gamma_k) \le 2$, and that the $H^1$-condition is satisfied. If $x_L = 0$ for all $x \in H^2(G)$,
does it follow that $L$ has a self-dual normal basis over $k$ ? This follows from Corollary 5.4 when $L$ is a field extension, in other words, when $\phi$ is surjective : indeed, then $H = G$.

\medskip
{\rm (ii)} The question above (see {\rm (i)}) has a negative answer for fields of higher cohomological dimensions. Indeed, by [BSe 94], III. 10.1,  there exist examples of $G$-Galois algebras $L$ over
fields of cohomological dimension 3 such that for all $n > 0$ we have $x_L = 0$ for all $x \in H^n (G)$, but $L$ does not have a self-dual normal basis
over $k$.
\medskip
{\rm (iii)} The converse of the question raised in  {\rm (i)} also has a negative answer : indeed, there exist examples of $G$-Galois algebras $L$ over ${\bf Q}$ having a self-dual normal basis such that there exists  $x \in H^2 (G)$ with  $x_L  \not = 0$ (see [BSe 94], III. 10.2).

\medskip
The following result was proved in  [BSe 94],
Corollaire 3.2.2  in the case where $k$ is an imaginary number field; the proof also applies for fields of cohomological dimension $\le 2$, using the results of [BP 95]. 
We give here an alternative proof.

\medskip
\noindent
{\bf Corollary 5.5.}  {\it Assume that ${\rm cd}_2(\Gamma_k) \le 2$, and that $$H^1(G) = H^2(G) = 0.$$
Then $L$ has a self-dual normal basis over $k$.}

\medskip
\noindent
{\bf Proof.} Since $H^1(G) = 0$, we have $G = G^2$. Let $A$ be orthogonal or unitary, and let us construct a central extension $V_A'$ of $G$ by ${\bf Z}/2{\bf Z}$, as follows. If $A$ is orthogonal, set $V'_A =V^G_A =  \tilde U_A(E_A) \times_{ U^0_{A}(E_A) }G$, with the notation of \S 3. If $A$ is unitary, then
we set $V'_A = V^G_A  =  F^{{\times}1}_{A} \times_{F^{{\times}1}_{A}} G$, the notation being as in \S4. In each case, we get a central extension $V_A'$ of $G$  by ${\bf Z}/2{\bf Z}$. Since
$H^2(G) = 0$, this extension is split. Note that the central extension $V_A$ of $H$ by ${\bf Z}/2{\bf Z}$ constructed in \S 3 and \S 4 is a subgroup of $V'_A$,
and that the restriction of the projection $V'_A \to G$ is the projection $V_A \to H$. Hence the extension $V_A$ is also split. This implies that we have 
$c_A(L) = 0$ for every orthogonal  $A$, and $d_A(L) = 0$ for every unitary $A$. By Theorem 5.3 this implies that $L$ has a self-dual
normal basis over $k$. 

\bigskip
\bigskip
{\bf \S 6. Cyclic groups of 2-power order}

\bigskip
In this section, $G$ is assumed to be  {\it cyclic  of order $2^n$}, with $n \ge 2$. We start by giving necessary and sufficient conditions for two $G$--Galois algebras to have isomorphic trace forms
in terms of cohomological invariants of degree 1 and 2 (see Proposition 6.2), namely the degree 1 invariants introduced in [BSe 94], and the discriminants of the
hermitian forms at the unitary factors (see \S 4). We then use the invariants defined
in the first part of \S 4 to give necessary and sufficient conditions for the existence of a self--dual normal basis.  We start with settling the case where $k$ contains the 4th roots
of unity :

\medskip
\noindent
{\bf Proposition 6.1.} {\it Assume that $k$ contains the $4$th roots of unity. Let $L$ and $L'$ be two $G$--Galois algebras. Then $q_L \simeq_G q_{L'}$ if and only if
$x_L = x_{L'}$ for all $x \in H^1(G)$.}

\medskip
\noindent
{\bf Proof.} The algebra $k[G]$ has two orthogonal factors $k$; since $k$ contains the 4th roots of unity, there are no other involution invariant factors. Therefore $u(L) = u(L')$
if and only if the cohomology classes $u$ associated to the two degree 1 orthogonal factors coincide, and this is equivalent with the condition $x_L = x_{L'}$ for all $x \in H^1(G)$.
Hence, by  [BSe 94], Proposition 1.5.1,  we have $q_L \simeq_G q_{L'}$.

\medskip

More generally, we have :
\medskip
\noindent
{\bf Proposition 6.2.} {\it Let $L$ and $L'$ be two $G$--Galois algebras. Then $q_L \simeq_G q_{L'}$ if and only if the following conditions hold :

\medskip

{\rm (i)} $x_L = x_{L'}$ for all $x \in H^1(G)$.

\medskip
{\rm (ii)}  ${\rm disc}_A(q_L) = {\rm disc}_A(q_{L'})$ for all unitary factors $A$ of $k[G]$.}

\bigskip From now on, we assume that $k$ {\it does not contain the $4$th roots of unity}. 
We start by introducing some notation. Set $A(i) = k[X]/(X^{2^{i-1}} + 1)$, for $i = 1,\dots,n$; then the factors of $k[G]$ are $k$, and $A(1),\dots,A(n)$. Note that
$k$ and $A(1)$ are orthogonal, and $A(2), \dots, A(n)$ are unitary. For
$i = 2,\dots,n$, we have $A(i) = F_{A(i)}$.


\medskip
\noindent
{\bf Proof of Proposition 6.2.}
For all factors $A$ of $k[G]$, let us denote by $h_A$, respectively $h'_A$, the hermitian form over $(A,\sigma_A)$ determined by $q_L$, respectively $q_{L'}$. 

\medskip
Assume that $q_L \simeq_G q_{L'}$. Then {\rm (i)} holds by [BSe 94], Proposition 2.2.1. Let $A$ be a unitary factor; then the hermitian forms $h_A$ and $h'_A$ are
isomorphic. Since ${\rm disc}_A(q_L)$ and $ {\rm disc}_A(q_{L'})$ are invariants of these hermitian forms, condition {\rm (ii)}  holds as well.

\medskip
Conversely, suppose that  {\rm (i)} and  {\rm (ii)} hold. Let us show that $u_A(L) = u_A(L')$ for all factors $A$. Condition  {\rm (i)} implies that this is true for
$A = k$ and $A = A(1)$; indeed, in both cases the group $U_A$ is of order 2. Let us assume that $A$ is a unitary factor, that is, $A = A(i)$ for some $i = 2,\dots,n$. Note that $A = F_A$, hence the hermitian forms $h_A$ and 
$h_A'$ are one dimensional hermitian forms over the commutative field $F_A$. Such a form is
determined up to isomorphism by its discriminant; hence condition {\rm (ii)} implies that $h_A \simeq h_A'$. Therefore we have $u_A(L) = u_A(L')$ for all factors $A$,
hence $u(L) = u(L')$,
and by [BSe 94], Proposition 1.5.1 we have $q_L \simeq_G q_{L'}$.

\medskip

Let us recall a notation from [Se 84], 1.5 or
[Se 92], 9.1.3 : we denote by $s_m \in H^2(S_m)$ the element of $H^2(S_m)$ corresponding to the central extension

$$1 \to {\bf Z}/2{\bf Z} \to \tilde S_m \to S_m \to 1$$ which is characterized by the properties :

\smallskip

1. A transposition in $S_m$ lifts to an element of order 2 in $\tilde S_m$.

2. A product of two disjoint transpositions lifts to an element of order 4 in $\tilde S_m$.

\medskip If $m$ is a power of 2, let us denote by  $C_m$ the cyclic group of order $m$, and by $e_m$ be the unique non-trivial element of $H^2(C_m)$. 
Sending a generator of $C_m$ to an $m$-cycle of $S_m$ defines an injective homomorphism $f : C_m \to S_m$; we denote by $f^* : H^2(S_m) \to H^2(C_m)$ the homomorphism
induced by $f$.  

\medskip
 If $q$ is a quadratic form over $k$, we denote by $w_2(q)$ its Hasse-Witt invariant (see
for instance [Se 84], 1.2 or [Se 92], 9.1.2); it is an element of $H^2(k)$.

\medskip
\noindent
{\bf Lemma 6.3.} {\it  Let $m$ be a power of 2. 

\medskip

{\rm (i)}  We have $f^*(s_m) = e_m$ in $H^2(C_m)$. 


\medskip
{\rm  (ii)} Let $\psi : \Gamma_k \to C_m$  be a continuous homomorphism, and let $K$ be the \'etale algebra over $k$ corresponding to $\phi$. Then the obstruction to the lifting of $\phi$   to a homomorphism 
$\Gamma_k \to C_{2m}$ is
$$w_2(q_K) + (2)(D_K)$$
where $D_K$ is the discriminant of $K$, and $(2)(D_K)$ denotes the cup product of the elements $(2)$ and $(D_K)$ of $H^1(k)$.}

\medskip
\noindent
{\bf Proof.} {\rm (i)}  Let $\tilde C_m$ be the inverse image of $C_m$ in $\tilde S_m$; it suffices to show that $\tilde C_m \simeq C_{2m}$, in other words that $\tilde C_m$ is a non-trivial extension of $C_m$. Raising an $m$-cycle of $S_m$ to the ${m \over 2}$-th power yields a product of ${m \over 2}$ disjoint transpositions, and the inverse image of such an element in
$\tilde S_m$ is of order 4. Hence $\tilde C_m$ is a non-trivial extension of $C_m$.

\medskip
{\rm (ii)} The obstruction to the lifting of $\psi$ is $\psi^*(e_m) \in H^2(k)$. Since $f^*(s_m) = e_m$
by {\rm (i)}, we have $$(f \circ \psi)^*(s_m) = \psi^*(e_m).$$  On the other hand,  $(f \circ \psi)^*(s_m) = w_2(q_K) + (2)(D_K)$ by [Se 84], Theorem 1.


\medskip
\noindent
{\bf Proposition 6.4.} {\it Let $L$ be a $G$-Galois algebra,  and assume that the $H^1$--condition holds. Then we have

\medskip
{\rm (i)} Let $A$ be a unitary factor of $k[G]$.  If $A \not = A(n)$, then $d_A(L) = 0$.

\medskip
{\rm (ii)} Let $L = K \times \cdots \times K$, where $K$ is a  field extension of $k$. Then  $$d_{A(n)}(L) = w_2(q_K) + (2)(D_K).$$}

\noindent
{\bf Proof.} Let $\phi : \Gamma_k \to G$ be a homomorphism associated to $L$, let $H = \phi (\Gamma_k)$, and let us denote by $|H|$ its order. Recall from \S 4 that the extension

$$(*) \ \ \ \ \ \ 1 \to {\bf Z}/2{\bf Z} \to V_A   \to H \to 1$$ is defined by $V_A = \{(x,h) \in F_A^{\times 1} \times H \ | \ x^2 = i_A(h) \}$. Let us show that this extension is split if $A \not = A(n)$.
Note that the group $V_A$ is abelian, and hence $(*)$ is not split if and only if $V_A$ is a cyclic group of order $2|H|$.  On the other hand, if $A \not = A(n)$, then the order of $i_A(H)$ is strictly less than $|H|$, hence the group $V_A$ does not have any elements of order $2|H|$. Therefore the extension $(*)$ is split, and hence $d_A(L) = 0$; this completes the proof of {\rm (i)}.

\medskip
If $A = A(n)$, then
the group $V_A$ is cyclic of order $2|H|$, and the extension $(*)$ is not split. Recall that we denote by $e_A \in H^2(H)$ the class of this
extension, and that $d_A = \phi^*(e_A) \in H^2(k)$. Note that $ \phi^*(e_A)$ is also the obstruction for the lifting of $\phi : \Gamma_k \to H$ to a 
continuous homomorphism $\Gamma_k \to V_A$; by Lemma 6.3 {\rm (ii)} this obstruction
is equal to $w_2(q_K) + (2)(D_K)$, hence {\rm (ii)}  is proved. 

\medskip
\noindent
{\bf Corollary 6.5.} {\it Let $L$ be a $G$-Galois algebra,  and assume that the $H^1$--condition holds. Then $L$ has a self--dual normal basis if and only if  
${\rm Res}_{E_{A(n)}/k}(d_{A(n)}(L)) = 0 \  {in} \  {\rm Br}_2(E_{A(n)}).$}

\medskip
\noindent
{\bf Proof.} Proposition 6.2 implies that $L$ has a self-dual normal basis if and only if the $H^1$-condition holds and if ${\rm disc}_A(q_L) = 0$ for all unitary factors
$A$ of $k[G]$. By Theorem 4.5 we have ${\rm Res}_{E_A/k}(d_A(L)) = {\rm disc}_A(q_L)$, and  Proposition 6.4 {\rm (i)}  implies that $d_A(L) = 0$ if $A \not = A(n)$. This completes
the proof of the corollary.

\medskip
\noindent
{\bf Corollary 6.6.} {\it Let $L$ be a $G$-Galois algebra,  and assume that the $H^1$--condition holds. 
Let $L = K \times \cdots \times K$, where $K$ is a  field extension of $k$, with $Gal(K/k)$  cyclic of order $m$. If $K$ can be embedded in a Galois extension
of $k$ with cyclic Galois group of order $2m$, then $L$ has a self-dual normal basis.}

\medskip
\noindent
{\bf Proof.} Assume that $K$ can be embedded  in a Galois extension
of $k$ with cyclic Galois group of order $2m$. Then by Lemma 6.3 {\rm (ii)} we have  $w_2(q_K) + (2)(D_K) = 0$.  By Proposition 6.4 {\rm (ii)}, this  implies that $d_{A(n)}(L) = 0$, and hence by Corollary 6.5 the $G$-Galois algebra $L$ has a self-dual normal basis. 

\medskip
\noindent
{\bf Example 6.7.} Assume that  $G$ is of order $8$.   
Let $a,b,c, \epsilon \in k$ with 
$a^2 - b^2 \epsilon = c^2 \epsilon$; assume $c$ non-zero, and  $\epsilon$ not a square. 
Set $x = \sqrt \epsilon$, and let $K = k(\sqrt {a + bx})$; note that $D_K = \epsilon$, and that  $K/k$ is a cyclic extension of degree $4$ (see for instance [Se 92], Theorem 1.2.1). 
Let $L$ be the $G$-Galois algebra induced from $K$. Let us prove that

\medskip

\centerline {$L$ has a self-dual normal basis $\iff$  $a$ is a sum of 4 squares in $k$. } 

\medskip
\noindent
Indeed, set $A = A(3)$; by 
Corollary 6.5
the $G$-Galois algebra $L$ has a self-dual normal basis if and only if ${\rm Res}_{E_{A}/k}(d_{A}(L)) = 0$. We have $d_A(L) = w_2(q_K) + (2)(\epsilon)$ by Proposition 6.4 {\rm (ii)}.

\medskip Let us show that $w_2(q_K) = (-1)(a)$. Set $y = \sqrt {a + bx}$. Then $\{1,x,y,xy\}$ is a basis of $K$ over $k$, and in this basis the quadratic form $q_K$ is the
orthogonal sum of the diagonal form $\langle 1,\epsilon \rangle$ and of the quadratic form $q$ given by $a X^2 + 2b \epsilon XY + a \epsilon Y^2$. The form  $q$ represents $a$,
and its determinant is $\epsilon(a^2 - b^2  \epsilon) = c^2 \epsilon^2$, hence ${\rm det}(q) = 1$ in $k^2/k^{\times 2}$. This implies that $q \simeq \langle a, a \rangle$,  hence
$q_K \simeq \langle 1,\epsilon,a,a \rangle$, and  $w_2(q_K) = (a)(a) = (-1)(a)$.

\medskip
Therefore $d_A(L) = (-1)(a) + (2)(\epsilon)$. Note that $E_A = k(\sqrt 2)$; hence ${\rm Res}_{A/k}(d_A(L)) = 
{\rm Res}_{k(\sqrt 2)/k}((-1)(a))$, and 
this element is $0$ if and only if $a$ is a sum of two squares in $k(\sqrt 2)$; or, equivalently, that $a$ is a sum of 4 squares in $k$ (see Example 5.2 {\rm (ii)}). 

\bigskip
Note that combining this example with Example 5.2 {\rm (i)} we get a necessary and sufficient condition for a $C_8$-Galois algebra to have a self-dual normal basis. 

\bigskip
\bigskip
{\bf \S 7. Self-dual normal bases over local fields}

\bigskip We keep the notation of the previous sections, and assume that $k$ is a (non-archimede\-an) local field. The aim of this section is to give a necessary and sufficient condition for the existence of self-dual normal
bases in terms of invariants defined over  $k$. 

\medskip
We say that $A$ is {\it split} if it is a matrix algebra over its center.

\medskip
\noindent
{\bf Theorem 7.1.} {\it The $G$-Galois algebra $L$ has a self-dual normal basis if and only if the $H^1$-condition holds, and 
\medskip
{\rm (i)} For all orthogonal $A$ such that $[E_A:k]$ is odd and $A$ is split, we have $c_A(L) = 0$  in ${\rm Br_2}(k)$. 
\medskip
{\rm (ii)} For all unitary  $A$  such that $[E_A:k]$ is odd, we have $d_A(L) = 0$ in ${\rm Br_2}(k)$. }

\medskip
\noindent
{\bf Proof.} Assume that the $H^1$-condition is satisfied and that {\rm (i)} and {\rm (ii)} hold. Note that if 
$A$ is not split, then we have
 ${\rm Br}_2(E_A)/\langle A \rangle = 0$, and that if 
$[E_A:k]$ is even, then the map ${\rm Res}_{E_A/k}  : {\rm Br_2}(k) \to {\rm Br_2}(E_A)$ is trivial. 
 Therefore for all orthogonal  $A$ we have
${\rm Res}_{E_A/k}(c_A(L)) = 0$
in ${\rm Br}_2(E_A)/\langle A \rangle$,
and for all  unitary 
$A$ we have ${\rm Res}_{E_A/k}(d_A(L)) = 0$ in  ${\rm Br}_2(E_A)$. By Theorem 5.3, this implies that $L$ has a self-dual normal basis.

\medskip Conversely, suppose that $L$ has a self-dual normal basis. Then the $H^1$-condition holds by Proposition 2.1. By Theorem 5.1 we have 
${\rm Res}_{E_A/k}(c_A(L)) = 0$
in ${\rm Br}_2(E_A)/\langle A \rangle$ for all orthogonal  $A$. Since ${\rm Res}_{E_A/k} : {\rm Br_2}(k) \to {\rm Br_2}(E_A)$ is injective if $[E_A:k]$ is odd,
condition {\rm (i)} holds. Moreover, Theorem 5.1  implies that if $A$ is unitary, then ${\rm Res}_{E_A/k}(d_A(L)) = 0$ in  ${\rm Br}_2(E_A)$. Applying
again the injectivity of ${\rm Res}_{E_A/k}$ when $[E_A:k]$ is odd, we obtain condition {\rm (ii)}. This completes the proof of the theorem. 

\bigskip
\bigskip
{\bf \S 8. Self-dual normal bases over global fields}

\bigskip  We keep the notation of the previous sections. Assume that $k$ is a global field, and let $\Omega_k$ be the set of places of $k$. For all $v \in \Omega_k$, we denote by $k_v$ the completion of $k$ at $v$. 
For all $k$-algebras $R$, set $R^v = R \otimes_k k_v$. We say that a $G$-Galois algebra is {\it split} if it is isomorphic to a direct product of copies of $k$ permuted
by $G$. We now apply the Hasse principle of [BPS 13] together with Theorem 7.1 above to give necessary and sufficient conditions for the existence of a self-dual
normal basis over $k$. 

\medskip
Note that the fields $E_A$ are abelian Galois extensions of $k$ (cf. 1.2).

\medskip
For all finite places $v$, let us write $E_A^v = K_A(v) \times \cdots  \times K_A(v)$, where $K_A(v)$ is a field extension of $k_v$. Set
$n^v_A = [K_A(v):k_v]$.

\medskip We need additional notation in the case when $A$ is unitary. 
Note that while $A$ is a central simple algebra over $F_A$, and $F_A/E_A$ is a quadratic
extension, for some places $v \in \Omega_k$ we may have $F^v_A = E^v_A\times E^v_A$  with $\sigma_A$ permuting the components, and $A^v = B  \times B$ for some $k_v$-algebra $B$. In order to take this
into account, we set $\epsilon^v_A = 0$ if $F^v_A = E^v_A\times E^v_A$, and $\epsilon^v_A = 1$ otherwise. 

\medskip
\noindent
{\bf Theorem 8.1.} {\it The $G$-Galois algebra $L$ has a self-dual normal basis if and only if the $H^1$-condition holds, if $L^v$ is split for all real places $v$, and
if for all finite places $v$ we have 
\medskip
{\rm (i)} For all orthogonal $A$ such that $n_{A}^v$ is odd  and $A^v$ is split, we have $c_A(L) = 0$  in ${\rm Br_2}(k_v)$. 
\medskip
{\rm (ii)} For all unitary  $A$  such that $n_{A}^v$ is odd and $\epsilon^v_A = 1$, we have $d_A(L) = 0$ in ${\rm Br_2}(k_v)$.}

\medskip
\noindent
{\bf Proof.} If $L$ has a self-dual normal basis, then $L^v$ is split for all real places $v$ by [BSe 94], Corollaire 3.1.2, and conditions {\rm (i)} and {\rm (ii)}  hold for
all finite places $v$ by Theorem 7.1. Conversely, assume that $L^v$ is split for all real places $v$, and that
for all finite places $v$ conditions {\rm (i)} and {\rm (ii)} hold. Then  [BSe 94], Corollary 3.1.2 (for real places) and Theorem 7.1 (for finite places) imply the existence of a self-dual normal
basis for $L^v$, for all $v \in \Omega_k$. By the Hasse principle result of [BPS 13], Theorem 1.3.1, the $G$-Galois algebra  $L$ has a self-dual normal basis over $k$.

\bigskip
\bigskip

\noindent
{\bf Bibliography} 

\bigskip \noindent 
[BL 90] E. Bayer-Fluckiger and H.W. Lenstra, Jr.,  Forms in odd degree
extensions and self-dual normal bases,
{\it Amer. J. Math.} {\bf 112} (1990), 359-373.

\bigskip \noindent
[B 94] E. Bayer-Fluckiger,
Multiplicateurs
de similitudes,
{\it C.R. Acad. Sci. Paris},
{\bf 319} (1994), 1151-1153.

\bigskip \noindent
[BF 15] E. Bayer-Fluckiger and U. First, Patching and weak approximation in isometry groups, {\it Trans. Amer. Math. Soc.}, to appear.


\bigskip \noindent 
[BP 95]  E. Bayer-Fluckiger and R. Parimala, Galois cohomology of the classical groups over fields of cohomological dimension $\le 2$,  {\it Invent. math.}  {\bf 122}  (1995),  195-229.
 

\bigskip \noindent
[BPS 13] E. Bayer-Fluckiger, R. Parimala and J-P. Serre, Hasse principle for $G$-trace forms, {\it Izvestiya Math.} {\bf 77} (2013), 437-460 (= {\it Izvestiya RAN}, {\bf 77} (2013), 5-28).

\bigskip \noindent
[BSe 94] E. Bayer-Fluckiger and J-P. Serre, Torsions quadratiques
et bases normales autoduales,
{\it Amer. J. Math.} {\bf 116} (1994),
1-64.

\bigskip \noindent
[K 69] M. Kneser, {\it Galois Cohomology of the Classical Groups}, Tata Lecture Notes in Mathematics {\bf 47} (1969). 

\bigskip \noindent
[KMRT 98] M. Knus, A. Merkurjev, M. Rost and  J-P. Tignol, {\it The Book of Involutions}, AMS Colloquium Publications {\bf 44}, 1998. 

\bigskip \noindent
[L 05] T. Y. Lam, {\it Introduction to quadratic forms over fields},  Graduate Studies in Mathematics {\bf 67} Amer. Math. Soc. Providence, RI (2005). 

\bigskip \noindent
[O 84] J. Oesterl\'e, Nombres de Tamagawa et groupes unipotents en caract\'eristique $p$, {\it Invent. math.} {\bf 78} (1984), 13-88.

\bigskip \noindent
[Sa 81] J-J. Sansuc, Groupe de Brauer et arithm\'etique des groupes algebriques lin\'eaires sur un corps de nombres, {\it J. reine angew. Math.} {\bf 327} (1981), 12-80.

\bigskip \noindent
[Sch 85] W. Scharlau, {\it Quadratic and Hermitian Forms}, Grundlehren der Math. Wiss.
Springer-Verlag (1985).

\bigskip\noindent 
[Se 64]  J-P. Serre, {\it Cohomologie galoisienne}, Lecture Notes in Mathematics, Springer-Verlag
(1964 and 1994).

\bigskip\noindent
[Se 68] J-P. Serre, {\it Corps locaux}, Hermann (1968).

\bigskip \noindent
[Se 84] J-P. Serre, L'invariant de Witt de la forme ${\rm Tr}(x^2)$, {\it Comment. Math. Helv.} {\bf 59} (1984), 651-676.

\bigskip \noindent
[Se 92] J-P. Serre, {\it Topics in Galois Theory}, Research Notes in Mathematics, Jones and Barlett Publishers (1992).

\bigskip
\bigskip
\bigskip

E. Bayer-Fluckiger

\'Ecole Polytechnique F\'ed\'erale de Lausanne

EPFL/FSB/MATHGEOM/CSAG

Station 8

1015 Lausanne, Switzerland

\medskip
eva.bayer@epfl.ch

\bigskip
\bigskip
R. Parimala

Department of Mathematics $ \&$ Computer Science

Emory University

Atlanta, GA 30322, USA.

\medskip
parimala@mathcs.emory.edu

\bye